\documentclass{amsart}
\usepackage{enumerate, amsfonts,
amssymb,xypic,amsmath,amsthm,pstricks,pst-node,pst-tree}

\newlength{\tabwidth}
\newlength{\tabheight}
\newlength{\tabrule}
\newlength{\tabwidthx}
\newlength{\tabheightx}

\def\gentabbox#1#2#3#4{\vbox to \tabheight{\setlength{\tabrule}{#3}%
  \setlength{\tabwidthx}{#1\tabwidth}\addtolength{\tabwidthx}{\tabrule}%

\setlength{\tabheightx}{#2\tabheight}\addtolength{\tabheightx}{-\tabheight}%
  \hbox to #1\tabwidth{%
 \hspace{-0.5\tabrule}\rule{\tabrule}{#2\tabheight}\hspace{-\tabrule}%
    \vbox to #2\tabheight{\hsize=\tabwidthx%
      \vspace{-0.5\tabrule}\hrule width\tabwidthx height\tabrule%
      \vspace{-0.5\tabrule}\vfil%
      \hbox to \tabwidthx{\hss#4\hss}%
        \vfil\vspace{-0.5\tabrule}%
      \hrule width\tabwidthx height\tabrule\vspace{-0.5\tabrule}}%
 \hspace{-\tabrule}\rule{\tabrule}{#2\tabheight}\hspace{-0.5\tabrule}}%
  \vspace{-\tabheightx}}}
\def\genblankbox#1#2{\vbox to \tabheight{\vfil\hbox to
#1\tabwidth{\hfil}}}
\def\tabbox#1#2#3{\gentabbox{#1}{#2}{0.4pt}{\strut #3}}

\catcode`\:=13 \catcode`\.=13 \catcode`\;=13
\catcode`\>=13 \catcode`\^=13
\def:#1\\{\hbox{$#1$}}
\def.#1{\tabbox{1}{1}{$#1$}}
\def>#1{\tabbox{2}{1}{$#1$}}
\def^#1{\tabbox{1}{2}{$#1$}}
\def;{\genblankbox{1}{1}\relax}
\catcode`\:=12 \catcode`\.=12 \catcode`\;=12
\catcode`\>=12 \catcode`\^=7

\newenvironment{tableau}{\bgroup\catcode`\:=13 \catcode`\.=13
  \catcode`\;=13 \catcode`\>=13 \catcode`\^=13
  \setlength{\tabheight}{3ex}\setlength{\tabwidth}{3ex}%
  \def\b##1##2##3{\gentabbox{##1}{##2}{1.2pt}{\vbox{##3}}}%
  \def\n##1##2##3{\gentabbox{##1}{##2}{0.4pt}{\vbox{##3}}}%
  \vbox\bgroup\offinterlineskip}{\egroup\egroup}

\newtheorem{theorem}{Theorem}[section]

\newtheorem{lemma}[theorem]{Lemma}
\newtheorem{proposition}[theorem]{Proposition}

\newtheorem*{theorem*}{}

\theoremstyle{definition}
\newtheorem{definition}[theorem]{Definition}

\theoremstyle{remark}
\newtheorem{claim}{Claim}[section]

\newtheorem{example}[theorem]{Example}

\begin{document}
\title{Equivalence Classes in the Weyl groups of type $B_n$}
\author{Thomas Pietraho}
\email{tpietrah@bowdoin.edu} \subjclass[2000]{20C08, 05E10}
\keywords{Unequal parameter Iwahori-Hecke algebra, Domino Tableaux,
Robinson-Schensted Algorithm}
\address{Department of Mathematics\\Bowdoin College\\Brunswick,
Maine 04011}
 \maketitle

\begin{abstract}
We consider two families of equivalence classes in the Weyl groups
of type $B_n$ which are suggested by the study of left cells in
unequal parameter Iwahori-Hecke algebras.  Both families are indexed
by a non-negative integer $r$.  It has been shown that the first
family coincides with left cells corresponding to the equal
parameter Iwahori-Hecke algebra when $r=0$; the equivalence classes
in the second family agree with left cells corresponding to a
special class of choices of unequal parameters when $r$ is
sufficiently large. Our main result shows that the two families of
equivalence classes coincide, suggesting the structure of left cells
for remaining choices of the Iwahori-Hecke algebra parameters.

\end{abstract}
\section{Introduction}

Consider a Weyl group $W$ of type $B_n$, that is, the
hyperoctahedral group $\mathbb{H}_n$. We will need two
generalizations of the Robinson-Schensted algorithm. Let $SDT_r(n)$
be the set of standard domino tableaux of size $n$ and rank $r$.
Then for every $r$, there is a map
$$G_r: W \rightarrow SDT_r(n) \times SDT_r(n)$$ which is a
bijection of $W$ with its image, the set of same-shape pairs of
standard domino tableaux (see \cite{garfinkle1} and
\cite{vanleeuwen:rank}). Subsequently, we will write $(S_r(w),
T_r(w))$ for the pair $G_r(w)$.  If we write $SBT(n)$ for the set of
standard bitableaux of size $n$, then it is possible to define
another algorithm
$$H: W \rightarrow SBT(n) \times SBT(n).$$
The image of $H$ is the set of same shape pairs of standard
bitableaux, and with this restriction, the map is again a bijection
\cite{stanton:white}. When $r \geq n-1$, there is a natural
identification between the sets $SDT_r(n)$ and $SBT(n)$, and the
corresponding bijections coincide.

As in \cite{lusztig:unequal}, a Coxeter group $W$ and a weight
function $L$ on $W$ can be used to partition $W$ into left cells.
When $W$ is of type $B_n$, the weight function can be identified
with a positive scalar $s$, and the description of the corresponding
left cells is known when $s=1$, i.e. the equal parameter case
\cite{garfinkle3}, when $s=\frac{1}{2}$ and $s=\frac{3}{2}$
\cite{bgil} and \cite{lusztig:leftcells}, and when $s>n-1$
\cite{bonnafe:iancu}. We examine two families of equivalence classes
on $W$ suggested by the above results, in an effort to interpolate
them for intermediate integer values of $s$.

We define the first equivalence relation on $W$ by letting $w \sim_r
y$ iff the right tableaux in the image of $G_r$, $T_r(w)$ and
$T_r(y)$, are related by moving through a set of open cycles.  The
second equivalence relation is defined using the right tableaux
obtained from both $G_r$ and $G_{r+1}$.  Let $w
\leftrightsquigarrow'_r y$ iff either $T_r(w)=T_r(y)$ or
$T_{r+1}(w)=T_{r+1}(y)$ and take  the equivalence relation
$\leftrightsquigarrow_r$ to be its transitive closure.  It is easy
to see that for sufficiently large $r$, the two relations are the
same. When $r \geq n-1$,  the set of open cycles defining $\sim_r$
is necessarily empty, implying $w \sim_r y$ iff $T_r(w)=T_r(y)$, and
further, that $w \sim_r y$ iff $w \leftrightsquigarrow_r y$. But
even more is true.

\vspace{.1in} \noindent {\bf Theorem 3.6} {\it   \hspace{.02in} Fix
a non-negative integer $r$.  For elements $w, y \in W,$ we have $w
\sim_r y$ iff $w \leftrightsquigarrow_r y$.} \vspace{.1in}

Both families of equivalence classes are derived from results on
left cells of unequal parameter Iwahori-Hecke algebras for certain
values of the parameter $s$.  The relation $\sim_0$, as well as the
notion of open cycles, was introduced by D.~Garfinkle
\cite{garfinkle1} and used to describe the left cells in the equal
parameter case, i.e. when the scalar $s=1$. For all values of $r
\geq n-1$, the equivalence classes of $\leftrightsquigarrow_r$ in
$W$ are the same, and correspond to the left cells in the unequal
parameter case when $s > n-1$ (see \cite{bonnafe:iancu}). This
result of Bonnaf\'{e} and Iancu originally described the left cells
using right bitableaux of $H$, but it is easy to reconcile with this
formulation by first, recalling the relationship between $H$ and
$G_r$ for large $r$ and second, noting that $T_r(w)=T_r(y)$ iff
$T_{r+1}(w)=T_{r+1}(y)$ for these values of $r$.

Since $\sim_0$ captures the left cell structure of $W$ when $s=1$
and $\leftrightsquigarrow_r$ describes the left cells when
$r=s-1>n-1$, the above result can be thought of as interpolating
these results and suggests that $\sim_r$ and
$\leftrightsquigarrow_r$ can be used to describe the left cell
structure for intermediate integral values of $s$ as well.  There is
one more piece of evidence in support of this.

Bonnaf\'{e}, Geck, Iancu, and Lam have conjectured that when $s$ is
non-integral, the left cells are determined by the right domino
tableaux for appropriate values of $r$ \cite{bgil}. Lusztig's
conjecture \cite{lusztig:unequal}(22.29) as well as the propositions
\cite{lusztig:unequal}(22.24) and \cite{lusztig:unequal}(22.25) then
suggest that the common refinement $\leftrightsquigarrow_{s-1}$ of
these cells describes the left cell structure for integral values of
$s$. The main result of this paper shows that this statement can be
rephrased using the somewhat simpler equivalence relations
$\sim_{s-1}$, and reconciles the conjecture with the original result
of Garfinkle, Barbasch and Vogan when $s=1$.

Some comments concerning the geometrical considerations related to
cells in this context are appropriate.  First of all, according to
Kazhdan-Lusztig theory, one would like to attach to each left cell a
distinguished Weyl group representation together with a finite
group. In the equal parameter case, this Weyl group representation
is special in the sense of Lusztig and the finite group is the
component group of the centralizer of an associated nilpotent orbit.
In the asymptotic case, all cells are irreducible and one can take
the corresponding finite groups to be trivial.

For the conjectured left cells in the intermediate cases, a choice
of a distinguished Weyl group representation, or in other words, a
choice of a distinguished partition among the shapes of tableaux of
elements in the conjectured two-sided cell (see
\cite{bgil}(Conjecture D)), is not clear.  In the equal parameter
case, the special representation attached to a left cell
$\mathcal{C}$ corresponds to a partition which appears as a tableaux
shape in all of the left cells contained in the two-sided cell of
$\mathcal{C}$.  Unfortunately, such a partition need not exist in
the unequal parameter case. For instance, when $r=2$, no single
partition appears among the tableaux shapes in every conjectural
left cell contained in the conjectural two-sided cell defined from
the partition $(4,3^2,1)$.

It is possible, however, to attach a finite group to each
conjectured left cell $\mathcal{C}$ in a manner reminiscent of the
equal parameter case.  The right tableaux of the involutions
$\mathcal{C} \cap \mathcal{C}^{-1}$ are all related by moving
through a subset of non-core open cycles. Hence $\mathcal{C} \cap
\mathcal{C}^{-1}$ can be endowed with the structure of an elementary
abelian group of order $2^c$, where $c$ is the number of non-core
open cycles in the right tableaux of elements of $\mathcal{C}$. In
the equal parameter case, McGovern has described the relationship of
$\mathcal{C} \cap \mathcal{C}^{-1}$ to a  quotient of the component
group mentioned above \cite{mcgovern:leftcells}.

Secondly, the map $G_0$, which classifies cells in the equal
parameter case, is related to a more geometric
Robinson-Schensted-type algorithm described in \cite{mcgovern:ssmap}
and \cite{pietraho:components} which classifies orbital varieties in
nilpotent orbits. The relationship between them is studied in
\cite{mcgovern:triangular}, establishing a bijection between cells
and the orbital varieties contained in special nilpotent orbits that
is well-behaved with respect to a certain partial order.  It would
be interesting to investigate whether a link exists between $G_r$
for $r \neq 0$ and this geometric algorithm; however, because
partitions of arbitrary rank do not naturally correspond to
nilpotent orbits, the existence of such a relationship is not
obvious.

Finally, Gordon and Martino have recently linked the combinatorics
of cells in the unequal parameter case with the geometry of the
Calogero-Moser space in \cite{gordon:calogero}.  There, the
nilpotent points of Calogero-Moser space are shown to correspond to
the combinatorially-defined conjectural two-sided cells.

\section{Preliminaries}

\subsection{Domino Tableaux}

A Young diagram $D$ is a finite left-justified array of squares
arranged with non-increasing row lengths.  A square in row $i$ and
column $j$ of the diagram will be denoted $s_{ij}$ so that $s_{11}$
is  the uppermost left square in the Young diagram below:

$$
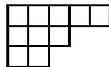
\begin{tiny}
\begin{tableau}
    :.{}.{}.{}.{}.{}\\
    :.{}.{}.{}\\
    :.{}.{}\\
\end{tableau}
\end{tiny}
$$

By  $\partial(D)$ of $D$, we will denote the set of squares $s_{ij}$
of $D$ such that either $s_{i,j+1}$ or $s_{i+1,j}$ does not lie in
$D$. We will also write $\rho(D)$ for the set of squares $s_{i,j+1}$
and $s_{i+1,j}$ that do not lie in  $D$ but $s_{ij} \in
\partial(D).$

Let $r \in \mathbb{N}$ and $\lambda$ be a partition of a positive
integer $m$; also write $\mathbb{N}_n=\{1, 2, \ldots, , n\}$.
     A {\it standard domino tableau of rank $r$ and shape $\lambda$} is a Young diagram
    of shape $\lambda$ whose squares are labeled by elements of $\mathbb{N}_n \cup \{0\}$ in such a
     way that the integer $0$ labels the square $s_{ij}$ iff $i+j<r+2$, each element of $\mathbb{N}_n$
      labels exactly two adjacent squares, and all labels increase
    weakly along both rows and columns.
We will write $SDT_r(\lambda)$ for the family of all domino tableaux
of rank $r$ and shape $\lambda$  and $SDT_r(n)$ for the family of
all domino tableaux of rank $r$ which contain exactly $n$ dominos.
The set of squares labeled by $0$ will be called the core of the
tableau.  We will write $\delta(T)$ for the set of $s_{ij}$ which
satisfy $i+j=r+2$,  and extend the notions of $\partial(D)$ and
$\rho(D)$ to tableaux, writing $\partial(T)$ and $\rho(T)$.

\subsection{Generalized Robinson-Schensted Algorithms}

We will consider the elements of the hyperoctahedral group $H_n$ as
subsets $w$ of $\mathbb{N}_n \times \mathbb{N}_n \times \{\pm 1\}$
with the property that the projections onto the first and second
components of $w$ are always bijections onto $\mathbb{N}_n$. The
element $w$ will be written as $\{(w_1,1,\epsilon_1), \ldots,
(w_n,n,\epsilon_n)\}$ and corresponds to the signed permutation
$(\epsilon_1 w_1, \epsilon_2 w_2, \ldots, \epsilon_n w_n)$.

We briefly describe the Robinson-Schensted bijections $G_r : H_n
\rightarrow SDT_r(n) \times SDT_r(n) ,$ following \cite{garfinkle1}
and \cite{vanleeuwen:rank}. The algorithm is based on a map $\alpha$
which inserts a domino with label $i$ into a domino tableau given an
element $(i,j, \epsilon)$ of $w \in H_n$. This insertion map is
similar to the usual Robinson-Schensted insertion map and is
precisely defined in \cite{garfinkle1}(1.2.5). To construct the left
tableau, start with $S_r(0)$, the only tableau in $SDT_r(0)$. Define
$S_r(1) = \alpha ((w_1,1,\epsilon_1),S_r(0))$ and continue
inductively by letting
$$S_r(k+1)=\alpha\big((w_{k+1},k+1,\epsilon_{k+1}),S_r(k)\big).$$
The left domino tableau $S_r(n)$ will be standard and of rank $r$.
The right tableau is defined to track the shape of the left tableau.
Begin by forming a domino tableau $T_r(1)$ by adding a domino with
label $1$ to $S_r(0)$ in such a way that $S_r(1)$ and $T_r(1)$ have
the same shape. Continue adding dominos by requiring that at each
step $T_r(k)$ lie in $SDT_r(k)$ and have the same shape as $S_r(k)$.
Again, the domino tableau $T_r(n)$ will be standard and of rank $r$.
We define the Robinson-Schensted map $G_r(w)=(S_r(n),T_r(n))$ and
will write $S_r(w)= S_r(n)$ and $T_r(w)=T_r(n)$ for the left and
right tableaux.

When $r=0$ or $1$,  the $G_r$ are precisely Garfinkle's algorithms;
for $r>1$ they are natural extensions to larger-rank tableaux. In
all cases, $G_r$ defines a bijection from $H_n$ to pairs of
same-shape tableaux in $SDT_r(n)$ (see \cite{vanleeuwen:rank}).
According to \cite{vanleeuwen:rank}, (4.2),  $G_r(w^{-1})=(T,S)$
whenever $G_r(w)=(S,T)$. In particular,  $w$ is an involution iff
$G_r(w)=(S,S)$ for some standard domino tableau $S$.

There is a natural description of the relationship between the
bijections $G_r$ for differing $r$ which we recount at the end of
the next section. We also point out that for $r$ sufficiently large,
$G_r$ recovers the algorithm $H$ mentioned in the introduction (see
\cite{stanton:white} and  \cite{okada}).

\subsection{Cycles}

We now review the notion of a cycle in a domino tableau.  It appears
in a number of references. See for instance \cite{garfinkle1} and
\cite{vanleeuwen:edge} as well as \cite{carre-leclerc} and
\cite{vanleeuwen:bijective}.

For $T \in SDT_r(n)$ we will call the  square $s_{ij}$ {\it fixed}
if $i+j$ has the opposite parity as $r$, otherwise, we'll call it
{\it variable}.  If $s_{ij}$ is variable and $i$ is odd, we will say
$s_{ij}$ is of type X; if $i$ is even, we will say $s_{ij}$ is of
type W. We will write $D(k,T)$ for the domino labeled by the
positive integer $k$ in $T$ and $supp \, D(k,T)$ will denote its
underlying squares. Write $label \, s_{ij}$ for the label of the
square $s_{ij}$ in $T$.  We extend this notion slightly by letting
$label \, s_{ij} =0$ if either $i$ or $j$ is less than or equal to
zero, and $label \, s_{ij} =\infty$ if $i$ and $j$ are positive but
$s_{ij}$ is not a square in $T$.

\begin{definition}
Suppose that  $supp \, D(k,T)= \{s_{ij},s_{i+1,j}\}$ or
$\{s_{i,j-1},s_{ij}\}$ and the square  $s_{ij}$ is fixed. Define
$D'(k)$ to be a domino labeled by the integer $k$ with $supp \,
D'(k,T)$ equal to
        $\{s_{ij}, s_{i-1,j}\}$     if $k< label \, s_{i-1,j+1}$ and
        $\{s_{ij}, s_{i,j+1}\}$    if $k> label \, s_{i-1,j+1}$.
Alternately, suppose that  $supp \, D(k,T)= \{s_{ij},s_{i-1,j}\}$ or
$\{s_{i,j+1},s_{ij}\}$ and the square  $s_{ij}$ is fixed. Define
$supp \, D'(k,T)$ to be
      $\{s_{ij},s_{i,j-1}\} $    if $k< label \, s_{i+1,j-1}$ and
         $\{s_{ij},s_{i+1,j}\}$         if $k> label \, s_{i+1,j-1}.$

\end{definition}

\begin{definition}
The cycle $c=c(k,T)$ through $k$ in a standard domino tableau $T$ is
a union of labels of dominos in $T$  defined by the condition that
$l \in c$ if either $l=k$, or either $supp \, D(l,T) \cap supp \,
D'(m,T) \neq
        \emptyset$ or $supp \, D'(l,T) \cap supp \,
D(m,T) \neq
        \emptyset$ for some $D(m,T) \in c$.
\end{definition}

We will often refer to the set of dominos with labels in a cycle $c$
as the cycle $c$ itself.   For a standard domino tableau $T$ of rank
$r$ and a cycle $c$ in $T$, define a domino tableau $MT(T,c)$ by
replacing every domino $D(l,T) \in c$ by the corresponding domino
$D'(l,T)$.

The tableau $MT(T,c)$ is standard, and in general, the shape of
$MT(T,c)$ will either equal the shape of $T$, or one square will be
removed (or added to the core) and one will be added
\cite{garfinkle1}, (1.5.27). A cycle $c$ is called closed in the
former case and open in the latter.   We will write $OC(T)$ for the
set of open cycles in $T$.  For $c\in OC(T)$, we will write $S_b(c)$
for the square that is either removed from the shape of $T$ or added
to the core of $T$ by moving through $c$. Similarly, we will write
$S_f(c)$ for the square that is added to the shape of $T$. Note that
$S_b(c)$ and $S_f(c)$ are always variable squares.

\begin{definition}
A variable square $s_{ij}$ in $\partial(T) \cup \rho(T)$ with the
property that neither $s_{i,j+1}$ nor $s_{i+1,j}$ lie in $T$ will be
called a {\it hole} if it is of type W and a {\it corner} if it is
of type $X$.
\end{definition}

 We will write $\Delta(T)$ for the
cycles through $\delta(T)$.  The squares $S_b(c)$ and $S_f(c)$ are
of the same type if $c \in \Delta(T)$.  However, for a cycle $c
\notin \Delta(T),$  one of $S_b(c)$ and $S_f(c)$ must be a corner
and the other a hole.  If the row number of $S_b(c)$ is smaller than
the row number of $S_f(c)$, we will call $c$ a down cycle;
otherwise, we will say $c$ is an up cycle.

Let $U$ be a set of cycles in $T$. According to \cite{garfinkle1},
(1.5.29), the order in which one moves through a set of cycles does
not matter, allowing us to unambiguously write $MT(T,U)$ for the
tableau obtained by moving-through all of the cycles in the set $U$.
Moving through a cycle in a pair of same-shape tableaux is a
slightly more delicate operation and requires the following
definition (see \cite{garfinkle2}(2.3.1)).

\begin{definition}
Consider $(S,T)$ a pair of same-shape domino tableaux, $k$ a label
in $S$, and $c$ the cycle in $S$ through $k$.  The extended cycle
$\tilde{c}$ of $k$ in $S$ relative to $T$ is a union of cycles in
$S$ which contains $c$.  Further, the union of two cycles $c_1 \cup
c_2$ lies in  $\tilde{c}$ if either is contained in $\tilde{c}$ and,
for some cycle $d$ in $T$,  $S_b(d)$ coincides with a square of
$c_1$ and $S_f(d)$ coincides with a square of $MT(S,c_2)$.  The
symmetric notion of an extended cycle in $T$ relative to $S$ is
defined in the natural way.
\end{definition}

 Let
$\tilde{c}$ be an extended cycle in $T$ relative to $S$. According
to the definition, it is possible to write $\tilde{c} =c_1 \cup
\ldots \cup c_m$ and find cycles $d_1, \ldots, d_m$ in $S$ such that
$S_b(c_i)=S_b(d_i)$ for all $i$, $S_f(d_m)=S_f(c_1)$, and
$S_f(d_i)=S_f(c_{i+1})$ for $1\leq i < m$. The union $\tilde{d}=d_1
\cup \dots \cup d_m$ is an extended cycle in $S$ relative to $T$
called the extended cycle corresponding to $\tilde{c}$.
Symmetrically, $\tilde{c}$ is the extended cycle corresponding to
$\tilde{d}$.

We define a moving through operation for a pair of same-shape domino
tableaux.  If we let $b$ be the ordered pair $(\tilde{c},\tilde{d})$
of extended cycles in $(S,T)$ that correspond to each other, then we
define $$MT((S,T), b)= (MT(S,\tilde{c}),MT(T,\tilde{d})).$$ As
desired, this operation produces another pair of same-shape domino
tableaux (\cite{garfinkle2}, (2.3.1)).

Finally, we are ready to describe the relationship between the
Robinson-Schensted maps $G_r$ introduced in the previous section.
Recall that the set $\delta(T)$ consists of the squares $s_{ij}$ of
$T$ adjacent to the core of $T$, and  $\Delta(T)$ is the set of open
cycles in $T$ that pass through the squares in $\delta(T)$. Let
$\gamma(S)$ be the extended cycles in $S$ relative to $T$ that pass
through $\delta(S)$; define $\gamma(T)$ similarly. If we write
$\gamma$ for the pair $(\gamma(S),\gamma(T))$, then
$$MMT((S,T))\equiv MT((S,T),\gamma)$$
is the minimal moving through map that clears all of the squares in
$\delta(S)$ and $\delta(T)$ simultaneously.

\begin{theorem} \cite{pietraho:rscore} Consider an element $w \in H_n$.  The Robinson-Schensted
  maps  $G_r$ and $G_{r+1}$ for rank $r$ and $r+1$ domino tableaux are related by
$$G_{r+1}(w) = MMT(G_{r}(w)).$$
\label{theorem:rscore}
\end{theorem}

\subsection{Cycle Structure}
\label{section:cyclestructure}

 Let us define a few objects that will be useful in describing the
 cycle structure of a domino tableau $T$.

\begin{definition}  Let $s_{ij}$ and $s_{kl}$ lie in the set of corners and holes of
$T$.  We will say that $s_{mn}$ is between $s_{ij}$ and $s_{kl}$ iff
$m$ is between $i$ and $k$ and $n$ is between $j$ and $l$ (where $m$
is between $i$ and $k$ iff $i \leq m \leq k$ or $i \geq m \geq k$).
We will also say $s_{ij}$ is above $s_{kl}$ if $i<k$.
\end{definition}

We will say that a cycle $d \in \Delta(T)$ is {\it adjacent} to a
cycle $c \in OC^*(T)=OC(T) \smallsetminus \Delta(T)$ if there is no
cycle $d' \in \Delta(T)$ such that $S_f(d')$ is between $S_f(d)$ and
$S_f(c)$.

The set of non-core open cycles $OC^*(T)$ has a partial order
defined by $c' \succeq c''$ iff $S_f(c'')$  is between $S_b(c')$ and
$S_f(c')$. Let $\mu(T) = \{c_1, c_2, \ldots\}$ be the set of maximal
elements in this poset.  For every $c \in \mu(T)$, we will write
$c_{\succeq}$ for the set of cycles smaller than or equal to $c$. We
also form a rooted tree $\tau(c)$ whose vertices correspond to
cycles $c' \in c_\succeq$, each labeled by $0$ if $c'$ is an up
cycle and $1$ if it is a down cycle. Edges in the tree $\tau(c)$ are
defined in the natural way from the Hasse diagram of the poset, and
for cycles of the same depth we place $c'$ to the left of $c''$ if
$S_b(c')$ is above $S_b(c'')$. Finally, $\tau(T)$ will denote the
ordered set of trees $(\tau(c_1), \tau(c_2), \ldots)$  where
$S_b(c_i)$ is above $S_b(c_{i+1}).$

\begin{example}
Consider the domino tableau
$$
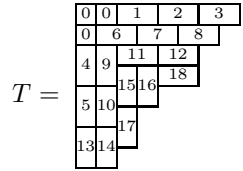
\raisebox{4ex}{$T=$ \;}
\begin{tiny}
\begin{tableau}
:.0.0>1>2>3\\
:.0>6>7>8\\
:^4^9>{11}>{12}\\
:;;^{15}^{16}>{18}\\
:^5^{10}\\
:;;^{17}\\
:^{13}^{14}\\
\end{tableau}
\end{tiny}
$$

\noindent It has three open cycles in the set $OC(T) \smallsetminus
\Delta(T)$: $c_1 = \{9, 10, 11,12,14\}$, $c_2=\{17\}$, and
$c_3=\{18\}$.  The set $\mu(T)$ contains only the cycle $c_1$ and
$\tau(T)=(\tau(c_1))$ where:

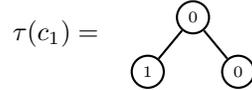
\begin{center}
\raisebox{-2ex}{$\tau(c_1)=$ \;}
\begin{tiny}
 \pstree[levelsep=8ex]{\Tcircle{0} }{
    \Tcircle{1}

    \Tcircle{0}
                    }
\end{tiny}
\end{center}

\end{example}

We will use $\tau(T)$ to keep track of the relative positions of the
non-core open cycles in $T$.   We will also need to keep track of
the exact locations of the beginning and final squares of all the
open cycles in $T$, for which we use the notion of a cycle structure
set.

\begin{definition} For a standard domino tableau $T$, we define the
cycle structure set of $T$ as the set of ordered pairs $cs(T)$
consisting of the beginning and final squares of every cycle in $T$.
That is:
$$cs(T) = \{ (S_b(c), S_f(c)) \; | \; c \in OC(T) \}.$$
\end{definition}

In general, the open cycles of a standard domino tableau can have
fairly complicated shapes.  However, for any $T$, it is always
possible to find another standard domino tableau with hook-shaped
open cycles and the same cycle structure set.

\begin{definition}  We will say that an open  cycle $c$ in $T$ is
{\it hook-shaped} iff the set of its underlying squares is entirely
contained in the union of one row and one column of $T$.
\end{definition}

\begin{proposition}\label{proposition:cyclestructure}
For every standard domino tableau $T$, there exists another standard
domino tableau $\Gamma(T)$ of the same shape such that:
    \begin{itemize}
        \item $cs(\Gamma(T)) = cs(T)$
        \item every $c \in OC(\Gamma(T))$ is hook-shaped.
    \end{itemize}
\end{proposition}

\begin{proof}
We first show that an appropriate domino tiling $d(T)$ of $shape(T)$
is possible and then show that its dominos can be labeled as a
standard tableau $\Gamma(T)$ with the required cycle structure.

We begin by assigning a hook in $shape(T)$ to every $c \in OC(T)$.
Write $S_b(c) = s_{ij}$ and $S_f(c) = s_{kl}$. For $c$ in
$\Delta(T)$, $k-i$ and $l-j$ are both even, implying that the hook
with ends $s_{ij}$ and $s_{k,l-1}$ and corner $s_{kj}$ as well as
the hook with ends $s_{ij}$ and $s_{k-1,l}$ and corner $s_{il}$ can
both be tiled by dominos. If the domino in the cycle $c$ of $T$
adjacent to the core is vertical, choose the former hook, otherwise,
choose the latter.  For $c$ in $OC^*(T)$, $k-i$ and $l-j$ are both
odd, again implying that the hook with ends $s_{ij}$ and $s_{k,l-1}$
and corner $s_{kj}$ as well as the hook with ends $s_{ij}$ and
$s_{k-1,l}$ and corner $s_{il}$ can be tiled by dominos.  If $c$ is
a down cycle, choose the former hook, otherwise, choose the latter.
We will write $h(c)$ for the set of dominos in the hook constructed
in this manner starting with a cycle $c \in OC(T).$  It is not hard
to see that it is possible to tile the above hook-shapes with
dominos for all cycles in $OC(T)$ simultaneously.

Next, we note that the remaining squares of $shape(T)$ can be tiled
with $2$-by-$2$ shapes.  To be more explicit, we make the following
definitions for the hooks in $shape(T)$ constructed above. Two such
hooks will be said to have adjacent vertical components if for every
pair of squares of the form $\{s_{ij}, s_{il}\}$ among the two
hooks, there is no other hook containing a square of the form
$s_{im}$ with $m$ between $j$ and $l$. A vertical component
containing a square $s_{ij}$ with a property that no other hook
contains a square of the form $s_{il}$ with $l$ greater than $j$
will be said to be adjacent to the boundary of $shape(T)$. We also
make the analogous definitions for the horizontal components of the
above hooks. Since adjacent corners and holes in $T$ have different
types, adjacent vertical components of two hooks must be separated
by an even number of columns, and vertical components adjacent to
the boundary of $shape(T)$ must be separated from it by an even
number of columns. Similarly, adjacent horizontal components of two
hooks must be separated by an even number of rows and horizontal
components adjacent to the boundary of $shape(T)$ must be separated
from it by an even number of rows.  This observation permits us to
tile the remaining squares of $shape(T)$ with $2$-by-$2$ shapes.

Using the above observations, define a domino tiling $d(T)$ of
$shape(T)$ as the set of $h(c)$ for all $c \in OC(T)$ together with
a tiling of the remaining $2$-by-$2$ shapes with pairs of adjacent
vertical dominos.

Next, we show that $d(T)$ can be numbered to create a standard
domino tableau.  To simplify the description, we first make the
following construction.

\begin{definition}
Let $e$ be a domino which intersects $\partial(D)$ of a domino-tiled
skew diagram $D$ of shape $\lambda \smallsetminus \mu$.  We
construct a rooted tree $t(e)$ with vertices corresponding to a
certain subset $v(e)$ of the dominos in the tiling of $D$. The
domino $e$ will correspond to the root of the tree.  The rest of
$t(e)$ is constructed recursively.

So suppose that $f$ is a domino of $D$ which corresponds to the
vertex $v$.  Whenever $f$ is either a horizontal domino containing
the corner $s_{ij}$ of a hook $h(c) $ with $c \in OC^*(T)$, or a
domino in a vertical portion of a hook $h(c)$ for $c \in OC(T)$ with
bottom square $s_{ij}$ fixed, we define $\tilde{f}$ as the domino
containing $s_{i+1,j-1}$.  We say that one domino is adjacent to
another if two of their underlying squares share a side.  The set of
children of $v$ in $t(e)$  consists of all the dominos adjacent to
the left side of $f$ together with the domino $\tilde{f}$ whenever
it is defined.  If there is only one child of $f$, place it as a
right child of $f$ in $t(e)$. Otherwise, order the children left to
right with the top-most domino as the left child. This process then
can be continued until the leaves of $t(e)$ correspond to dominos
adjacent to the boundary of $\mu$.
\end{definition}

We are now ready to describe how to number the dominos of $d(T)$.
Consider the domino $e$ at the top edge of $\partial(T)$.  The
vertices of $t(e)$ can be numbered with the integers $\{1, 2,
\ldots, |t(e)|\}$ according to a postfix order traversal, and the
labels of $t(e)$ can be transferred to $d(T)$. The entire process
can then be repeated for the skew diagram $d(T) \setminus v(e)$ with
labels starting at $|t(e)|+1$. Iterating this procedure until $d(T)$
is exhausted yields a numbered domino tiling of $d(T)$.  It is not
difficult to verify that it is in fact a standard domino tableau
which we call $\Gamma(T)$.

Finally, we check that $\Gamma(T)$ has the desired cycles. First, we
show that all dominos of $\Gamma(T)$ not contained in a hook $h(c)$
for some $c \in OC(T)$ lie in a closed cycle of size $2$ in
$\Gamma(T)$. Consider the dominos of $\Omega= \Gamma(T)
\smallsetminus \{h(c) \; | \; c \in OC(T)\}$. Choose a domino $e_1
\in \Omega$ whose top and left
 edges are not adjacent to any other members of $\Omega$. Its top
square $s_{ij}$ is necessarily fixed.  Further, the label of
$s_{i-1,j+1}$ is necessarily smaller than the label of $s_{ij}$
since it either lies outside of $shape(T)$, or the domino containing
it corresponds to either a vertex on a prior branch of its tree, or
even a tree that has been labeled previously. Hence $MT(e_1,
\Gamma(T))$ consists of the squares $\{s_{ij},s_{i,j+1}\}$. The
domino $e_2= \{s_{i,j+1},s_{i+1,j+1}\}$ of $\Gamma(T)$ has fixed
square $s_{i+1,j+1}$. Its label is necessarily smaller then the
label of $s_{i+2,j}$, since the latter square either lies outside
the Young diagram underlying $shape(T)$, or the domino containing it
corresponds to either a vertex on a later branch of its tree, or
even a tree labeled later. Hence $MT(e_2, \Gamma(T))$ consists of
the squares $\{s_{i+1,j},s_{i+1,j+1}\}$ and $\{e_1, e_2\}$ is a
closed cycle in $\Gamma(T)$.  This argument can be repeated again
with $\Omega$ replaced with $\Omega \smallsetminus \{e_1, e_2\}$
until $\Omega$ is exhausted.

Now, consider a hook $h(c)$ corresponding to $c \in OC(T)$.  When
labeled, it constitutes an open cycle of $\Gamma(T)$;  the fact that
if $e$ is a domino in $h(c)$, then the squares of $MT(e,T)$ again
lie in $h(c)$ follows directly from our method of numbering $d(T)$.
Furthermore, if $c \in OC(T)$, then $S_b(h(c)) = s_{ij} =S_b(c)$ and
$S_f(h(c)) = s_{kl}=S_f(c).$ Since the $h(c)$ are the only open
cycles of $\Gamma(T)$, we can conclude that $cs(\Gamma(T)) = cs(T)$.
\end{proof}

\begin{example}  We provide an example of the method of labeling the dominos
in a  domino diagram to form a standard domino tableau used in the
above construction. Consider the domino diagram $E$ of shape
$\{5,5,4,3\} \setminus \{2,1\}$. The domino $e$ is at the top of its
boundary:

$$
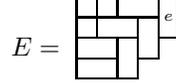
\raisebox{2ex}{$E=$ \;}
\begin{tiny}
\begin{tableau}
:.{}.{}>{}^{e}\\
:.{}>{}^{}\\
:>{}^{}\\
:>{}\\
\end{tableau}
\end{tiny}
$$

\noindent The dominos in $E$ can be represented in the rooted tree
$t(e)$ which can be numbered according to a postfix traversal order
as below.

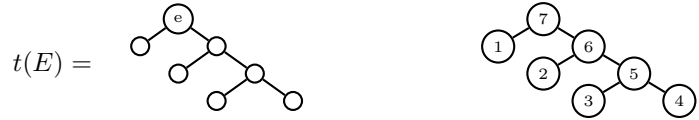
\begin{center}
\raisebox{-4ex}{$t(E)=$ \;}
\begin{tiny}
\raisebox{1ex}{ \pstree[levelsep=4ex]{\Tcircle{e} }{
    \Tcircle{}

    \pstree[levelsep=4ex]{\Tcircle{} }{
    \Tcircle{}
    \pstree[levelsep=4ex]{\Tcircle{} }{
          \Tcircle{}
          \Tcircle{}          }}}

}
\end{tiny}
\hspace{.7in}
\begin{tiny}
\raisebox{1ex}{ \pstree[levelsep=4ex]{\Tcircle{7} }{
    \Tcircle{1}

    \pstree[levelsep=4ex]{\Tcircle{6} }{
    \Tcircle{2}
    \pstree[levelsep=4ex]{\Tcircle{5} }{
          \Tcircle{3}
          \Tcircle{4}          }}}

}
\end{tiny}
\end{center}

\noindent Finally, this numbering of the vertices of the above tree
yields the standard domino tableau $T$:

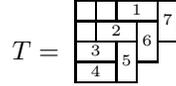
$$
\raisebox{2ex}{$T=$ \;}
\begin{tiny}
\begin{tableau}
:.{}.{}>{1}^{7}\\
:.{}>{2}^{6}\\
:>{3}^{5}\\
:>{4}\\
\end{tableau}
\end{tiny}
$$

\end{example}

Note that in the proof of the previous proposition, the only data
required from $T$ to construct $\Gamma(T)$ was its shape and the
cycle structure set $cs(T)$. Thus given a shape $\lambda$ of a
domino tableau and a set $cs$ of pairs of beginning and final
squares, we will write $\Gamma(\lambda, cs)$ for the standard domino
tableau constructed via the above process.

\section{Equivalence Relations on W}

We define two equivalence relations on domino tableaux in
$SDT_r(n)$.  Both can be used to define equivalence classes in $W$
of type $B_n$ by using the generalized Robinson-Schensted algorithms
$G_r$. The main result is that the equivalence classes thus defined
are in fact the same.

\begin{definition}\label{definition:sim}  Consider $T, T' \in SDT_r(n)$.  We say $T\sim_r T'$
iff there is a subset of non-core open cycles $U \subset OC^*(T')$
such that $T=MT(T',U)$.
\end{definition}

It is not difficult to see that $\sim_r$ defines an equivalence
relation on $SDT_r(n)$.  We will use it to define an equivalence
relation on $W$.  First, note that elements of an
$\sim_r$-equivalence class of $T$ correspond to $\{0,1\}$-labelings
of the vertices of the trees underlying $\tau(T)$. Such a labeling
determines a unique standard domino tableau $T'$ in the equivalence
class of $T$.

\begin{definition}\label{definition:equivalence1}  Consider $w, y \in W$ and let
 $w \sim_r y$ iff
$T_r(w) \sim_r T_r(y)$.
\end{definition}

Hence the equivalence relation $\sim_r$ on $W$ is completely
determined by an equivalence relation on right tableaux.

\begin{definition}\label{definition:leftrightsquigarrow}  Consider $T, T' \in SDT_r(n)$.  We will say $T \leftrightsquigarrow'_r T'$
iff there exist $S, S'\in SDT_r(n)$ such that $MMT(S,T)$ and
$MMT(S',T')$ have the same right tableau. We will let
$\leftrightsquigarrow_r$ be its transitive closure.
\end{definition}

\begin{definition}
 Consider $w, y \in W$ and define a relation $w \leftrightsquigarrow'_r y$ iff
$T_r(w) = T_r(y)$ or $T_{r+1}(w) = T_{r+1}(y)$.  Let
$\leftrightsquigarrow_r$ be its transitive closure.
\end{definition}

Both of the relations denoted by $\leftrightsquigarrow_r$ are in
fact equivalence relations.  Again, the equivalence relation
$\leftrightsquigarrow_r$ on $W$ can be completely expressed as an
equivalence relation on right tableaux, mirroring Definition
\ref{definition:equivalence1}.

\begin{proposition}\label{proposition:equivalence}  Consider $w, y \in W,$
then $w \leftrightsquigarrow_r y$
iff $T_r(w) \leftrightsquigarrow_r T_r(y)$.
\end{proposition}

\begin{proof}
Let's simplify the notation and write $\sim$ for $\sim_r$ and
$\leftrightsquigarrow$ for $\leftrightsquigarrow_r$. Note that $w
\leftrightsquigarrow' y$ means that either $T_r(w) = T_r(y)$, in
which case $T_r(w) \leftrightsquigarrow' T_r(y)$, or that
$T_{r+1}(w) = T_{r+1}(y)$. In the latter case, the description of
the map $MMT$ implies that $MMT(G_r(w))$ $MMT(G_r(y))$ have the same
right tableaux and $T_r(w) \leftrightsquigarrow' T_r(y)$. Hence the
relation $\leftrightsquigarrow$ on $W$ implies the relation
$\leftrightsquigarrow$ on tableaux.

Conversely, suppose that $T_r(w) \leftrightsquigarrow' T_r(y)$.
Consider tableaux $S$ and $S'$ such that $MMT(S, T_r(w))$ and
$MMT(S',T_r(y)$ have the same right tableaux.  If we let
$w'=G_r^{-1}(S,T_r(w))$ and $y'=G_r^{-1}(S',T_r(y))$, then
$w'\leftrightsquigarrow' y'$ since $G_{r+1}(w')= MMT(S, T_r(w))$ and
$G_{r+1}(y')= MMT(S', T_r(y)).$ Finally, $T_r(w)=T_r(w')$ and
$T_r(y)=T_r(y')$, implying $w \leftrightsquigarrow' w'
\leftrightsquigarrow' y'\leftrightsquigarrow' y$. Thus, the
equivalence relation $\leftrightsquigarrow$ on tableaux implies the
relation $\leftrightsquigarrow$ on $W$.
\end{proof}

Our main result states that $\sim_r$ and $\leftrightsquigarrow_r$ in
fact define the same equivalence classes in $W$.

\begin{theorem}\label{theorem:main}
Consider $w, y \in W$, then $w \sim_r y$ iff $w
\leftrightsquigarrow_r y$.
\end{theorem}
\begin{proof}
Again, let's  write $\sim$ for $\sim_r$ and $\leftrightsquigarrow$
for $\leftrightsquigarrow_r$. The fact that $w\leftrightsquigarrow
y$ implies $w \sim y$ is a consequence of the relationship between
the Robinson-Schensted maps $G_r$ and $G_{r+1}$ described in Theorem
\ref{theorem:rscore}.

If $w\leftrightsquigarrow ' y$, then either $T_r(w)=T_r(y)$ or
$T_{r+1}(w)=T_{r+1}(y)$.  In the former case, we automatically have
 $w \sim y$.  In the latter case,
$T_{r+1}(w)=T_{r+1}(y)$ and Theorem \ref{theorem:rscore} implies
that $MMT(G_r(w))$ and $MMT(G_r(y))$ have the same right tableaux.
According to the definition of the $MMT$ map,
$$T_{r+1}(w)= MT(T_r(w), \Delta(T_r(w)) \cup X)$$ for a subset of
non-core open cycles $X \subset OC^*(T_r(w))$ and
$$T_{r+1}(y) = MT(T_r(y), \Delta(T_r(y)) \cup Y)$$ for a subset of
non-core open cycles $Y \subset OC^*(T_r(y)).$ Since the
moving-through operation can be performed on disjoint sets of cycles
independently, $MT(T_r(w), X)$ must equal $MT(T_r(y),Y)$ . In
particular, this forces $T_r(w) = MT(T_r(y), U)$ for some subset $U
\subset OC^*(T_r(y))$, implying $T_r(w) \sim T_r(y)$ and $w \sim y$.

Our proof that $w \sim y$ implies $w \leftrightsquigarrow y$
requires the following lemma, whose proof we relate in the following
section.

\begin{lemma}\label{lemma:main}
Suppose that $T'' = MT(T', U)$ for some subset $U \subset OC^*(T')$,
and that the labels of $\tau(T')$ and $\tau(T'')$ disagree only on
one tree.  Then $T' \leftrightsquigarrow T''$.
\end{lemma}

Let $T = T_r(w)$ and $T'=T_r(y)$.  By definition, $w \sim y$ implies
$T' = MT(T, U)$ for some subset $U$ of non-core open cycles of $T$.
Let's suppose that the positions of the cycles of these tableaux are
given by the trees $\tau(T)=(\tau_1, \tau_2, \ldots, \tau_m)$ and
$\tau(T')=(\tau'_1, \tau'_2, \ldots, \tau'_m)$. Using the note
following Definition \ref{definition:sim}, construct domino tableaux
$T_i$ with $\tau(T_i) = (\tau_1, \tau_2, \ldots, \tau_i,
\tau'_{i+1}, \ldots \tau'_m)$. According to Lemma \ref{lemma:main},
we have constructed a sequence of tableaux $\{T_i\} $ satisfying
$$T' =T_0 \leftrightsquigarrow T_1 \leftrightsquigarrow T_2
\leftrightsquigarrow \ldots \leftrightsquigarrow T_m = T.$$ The
characterization of $\leftrightsquigarrow$ given in Proposition
\ref{proposition:equivalence} implies $w \leftrightsquigarrow y$, as
desired.
\end{proof}

\section{Main Lemma}
This section details the proof of Lemma \ref{lemma:main}.  We retain
the notation used therein and begin with a definition and a result
on $\{0,1\}$-labeled trees.

\begin{definition}
Suppose $\tau'$ and $\tau''$ are \{0,1\}-labeled trees as in Section
\ref{section:cyclestructure} which share the same underlying
unlabeled embedded rooted tree $\varsigma$.  For $\epsilon \in
\{0,1\}$, we say $\tau' \leftrightsquigarrow'_\epsilon \tau''$ if
there is a path $\kappa$ in $\varsigma$ satisfying
    \begin{enumerate}
        \item $\kappa$ contains the root of $\varsigma$
        \item labels of $\tau'$ and $\tau''$ agree on the vertices
        of $\varsigma \smallsetminus\kappa$,
        \item except for the vertex of greatest depth in $\kappa$, the labels of
        vertices in  $\kappa$
        in $\tau'$ and $\tau''$ agree, alternate, and begin with
        $\epsilon$ at the root.
   \end{enumerate}
We will write $\leftrightsquigarrow_\epsilon$ for the transitive
closure of $\leftrightsquigarrow'_\epsilon$.
\end{definition}

We note that according to the above, a path $\kappa$ can contain
just the root of $\varsigma$.  The next result shows that for each
$\varsigma$, there is only one equivalence class of
$\{0,1\}$-labeled trees.  While this means that the relation
$\leftrightsquigarrow_\epsilon$ is trivial, we've introduced it as
we will need to know the sequence of
$\leftrightsquigarrow'_\epsilon$ that accomplishes this.
\begin{proposition}\label{proposition:tree}
Two $\{0,1\}$-labeled  trees which share the same underlying rooted
tree are $\leftrightsquigarrow_\epsilon$-equivalent.
\end{proposition}

\begin{proof}
Consider two $\{0,1\}$-labeled trees $\tau'$ and $\tau''$ which
share the same underlying unlabeled rooted tree $\varsigma$. Let
$\tau_a$ be a labeling of the vertices of $\varsigma$ with vertices
of even depth labeled by $\epsilon$ and vertices of odd depth
labeled by $1-\epsilon$. Starting with the set of vertices of
maximal depth $m$ in $\varsigma$ and using the definition of
$\leftrightsquigarrow'_\epsilon$, we can find a sequence of trees
which agree with $\tau_a$ on vertices of depth less than $m$ and are
$\leftrightsquigarrow_\epsilon$-equivalent to the tree whose labels
agree with $\tau_a$ on vertices of depth less than $m$ and vertices
of $\tau'$ on maximal depth vertices.  This procedure can be
repeated successively for smaller depths, creating a sequence
$$\tau' = \tau_1 \leftrightsquigarrow'_\epsilon \tau_2
\leftrightsquigarrow'_\epsilon \ldots \leftrightsquigarrow'_\epsilon
\tau_l = \tau_a.$$ The above procedure can be repeated with $\tau'$
replaced with $\tau''$, finally creating the desired sequence of
labeled embedded trees
$$\tau' = \tau_1 \leftrightsquigarrow'_\epsilon \tau_2
\leftrightsquigarrow'_\epsilon \ldots \leftrightsquigarrow'_\epsilon
\tau_m = \tau''.$$

\end{proof}

\begin{example}
We exhibit the above procedure for the trees $\tau'$ and $\tau''$
and $\epsilon=1$.  Here, $\tau_a$ is the third tree in the sequence.
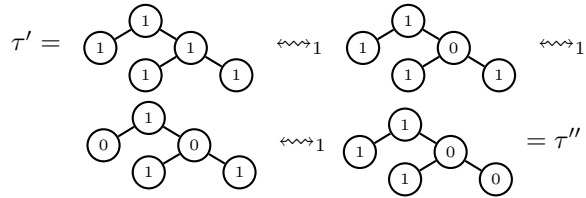
\begin{center}
\raisebox{-2ex}{$\tau'  =$}
\begin{tiny}
\raisebox{1ex}{
\pstree[levelsep=4ex]{\Tcircle{1} }{
    \Tcircle{1}

    \pstree[levelsep=4ex]{\Tcircle{1} }{
          \Tcircle{1}
          \Tcircle{1}          }}

}
\end{tiny}
\raisebox{-2ex}{$\leftrightsquigarrow_1$}
\begin{tiny}
\raisebox{1ex}{ \pstree[levelsep=4ex]{\Tcircle{1} }{
    \Tcircle{1}

    \pstree[levelsep=4ex]{\Tcircle{0} }{
          \Tcircle{1}
          \Tcircle{1}          }}

} \end{tiny} \raisebox{-2ex}{$\leftrightsquigarrow_1$}
\end{center}

\begin{center}
\raisebox{-2ex}{\phantom{$\leftrightsquigarrow_1  $}}
\begin{tiny}
\raisebox{1ex}{ \pstree[levelsep=4ex]{\Tcircle{1} }{
    \Tcircle{0}

    \pstree[levelsep=4ex]{\Tcircle{0} }{
          \Tcircle{1}
          \Tcircle{1}          }}

} \end{tiny} \raisebox{-2ex}{$\leftrightsquigarrow_1$}
\begin{tiny}
{\pstree[levelsep=4ex]{\Tcircle{1} }{
    \Tcircle{1}

    \pstree[levelsep=4ex]{\Tcircle{0} }{
          \Tcircle{1}
          \Tcircle{0}          }}
}
\end{tiny}\raisebox{-2ex}{$=\tau''$}
\end{center}
\end{example}

We are ready to start the  proof of the main lemma.   Suppose that
$T'' = MT(T', U)$ for some subset $U \subset OC^*(T')$, and that the
labels of $\tau(T')$ and $\tau(T'')$ disagree only on the $(k+1)$-st
tree:

$$\tau(T') = (\tau^1, \tau^2, \ldots, \tau^k, \tau', \tau^{k+2}, \ldots)$$
$$\tau(T'') = (\tau^1, \tau^2, \ldots, \tau^k, \tau'', \tau^{k+2}, \ldots)$$

Write $c$ for the cycle of $T'$ and $T''$ that corresponds to the
root of both the trees $\tau'$ and $\tau''$, and choose a cycle $d
\in \Delta(T')= \Delta(T'')$ adjacent to $c$ (see Section
\ref{section:cyclestructure}).  If $S_f(d)$ is above
 $S_f(c)$, then take $\epsilon = 1$. Otherwise, let $\epsilon=0$.
The previous proposition gives $\tau' \leftrightsquigarrow_\epsilon
\tau''$ via  the sequence
$$\tau' = \tau_1 \leftrightsquigarrow'_\epsilon \tau_2
\leftrightsquigarrow'_\epsilon \ldots \leftrightsquigarrow'_\epsilon
\tau_m = \tau''.$$ Using the note following Definition
\ref{definition:sim}, we can construct a sequence of tableaux
$\{T_i\}_{i=1}^m$ from this sequence of labeled embedded trees, each
satisfying $$\tau(T_i)=(\tau^1, \tau^2, \ldots, \tau^k, \tau_i,
\tau^{k+2}, \ldots).$$

\begin{claim}\label{claim}
With $\epsilon$ as above,  $T_i \leftrightsquigarrow' T_{i+1}$.
\end{claim}

\begin{proof}

Since $\tau_i \leftrightsquigarrow'_\epsilon \tau_{i+1}$, we must
have $T_{i+1} = MT(T_i, \tilde{c})$ for some cycle $\tilde{c}
\preceq c$.  Let $c=c_1, c_2, \ldots c_l=\tilde{c}$ be the maximal
chain between them in the poset of open cycles. Except perhaps for
$\tilde{c}$, it alternates between up and down cycles in both $T_i$
and $T_{i+1}$ because of the way $\leftrightsquigarrow'$ is defined.
Furthermore, unless $l=1$, $c_1$ must be an up cycle if $\epsilon=1$
and a down cycle if $\epsilon=0$. Let $\mathcal{B}$ be the cycle
structure set of $T_i$ with pairs of squares that correspond to
cycles in $\tau_i$ deleted. Write $c_0 =d$ and for $k \geq 1$, let
$$\mathcal{C}_k = (S_b(c_0),S_f(c_1))   \cup \bigcup_{1 \leq j \leq k/2}
(S_b(c_{2j}),S_f(c_{2j-2})) \cup \bigcup_{1 \leq j < k/2}
(S_b(c_{2j-1}), S_f(c_{2j+1}))$$ and
$$\mathcal{D}_k = \left\{ \begin{array}{ll}
\mathcal{B} \cup \mathcal{C}_k \cup (S_b(c_k), S_f(c_{k-1})) & \textrm{if $k$ is odd}\\
\mathcal{B} \cup \mathcal{C}_k \cup (S_b(c_{k-1}), S_f(c_{k}))& \textrm{if $k$ is even.}\\
\end{array} \right.
$$
When $k=0$, let $\mathcal{D}_0 = \mathcal{B}$. There are two
possibilities for the position of $c_l$ in $T_i$. First assume that
$c_1, c_2, \ldots, c_l$ alternates between up and down cycles.
Following Proposition \ref{proposition:cyclestructure} and
especially the comment after it, we can find standard domino
tableaux $S_i=\Gamma(shape(T_i),\mathcal{D}_{l})$ and
$S'_i=\Gamma(shape(T_{i+1}),\mathcal{D}_{l-1})$ satisfying $cs(S_i)
= \mathcal{D}_{l}$ and $cs(S'_i) = \mathcal{D}_{l-1}.$ By definition
of extended cycles and our choice of $\epsilon$,
\begin{align*}
ec(d,T_i,S_i) & = d \cup c_1 \cup c_2 \cup \ldots c_l, \text{ and }\\
ec(d,T_{i+1},S'_i) & =d \cup c_1 \cup c_2 \cup \ldots c_{l-1}.
\end{align*}
Hence $MMT(S_i,T_i)$ will have right tableau  $MT(T_i, \Delta(T_i)
\cup c_1 \cup c_2 \cup \ldots c_l)$ while $MMT(S'_i,T_{i+1})$ will
have right tableau  $MT(T_{i+1}, \Delta(T_{i+1}) \cup c_1 \cup c_2
\cup \ldots c_{l-1})$.  Since $\Delta(T_{i+1})=\Delta(T_{i})$ and
$T_{i+1} = MT(T_i, c_l)$, Definition
\ref{definition:leftrightsquigarrow} implies that $T_i
\leftrightsquigarrow' T_{i+1},$  as desired.

The second possibility for the position of $c_l$ in $T_i$ is that
both $c_{l-1}$ and $c_{l}$ are either up or down cycles.  The proof
follows as above, but with $S_i$ and $S'_i$ defined to satisfy
$cs(S_i) = \mathcal{D}_{l-1}$ and $cs(S'_i) = \mathcal{D}_{l}.$

\end{proof}

Armed with the above claim, we have found a sequence of tableaux
$\{T_i\}_{i=1}^m$ such that $T_1=T'$, $T_m=T''$ and $T_i
\leftrightsquigarrow' T_{i+1}$ for all $i$, verifying the main
lemma.

\begin{example}
We conclude with an example which will hopefully clarify the above
procedure.  Let
$$
\raisebox{3ex}{$T'=$ \;}
\begin{tiny}
\begin{tableau}
:.0>1>2\\
:>3>4\\
:>5>6\\
:^7>8\\
\end{tableau}
\end{tiny}
\raisebox{3ex}{\hspace{.3in} and \hspace{.3in} $T''=$ \;}
\begin{tiny}
\begin{tableau}
:.0>1>2\\
:>3>4\\
:^5>6\\
:;^8\\
:^7\\
\end{tableau}
\end{tiny}
$$
so that $T''=MT(T',U),$ where $U$ consists of all non-core open
cycles of $T'$.   We would like to show that $T'
\leftrightsquigarrow T''$, and so we need to construct a sequence of
$\leftrightsquigarrow'$-equivalent tableaux. Both $\tau(T')$ and
$\tau(T'')$ contain only one tree, so we are in the setting of the
Lemma \ref{lemma:main}.   We first find a sequence of labeled trees
which exhibit $\tau(T') \leftrightsquigarrow_1 \tau(T'')$.  The
proof of Proposition \ref{proposition:tree} yields the sequence
\begin{center}
\raisebox{-2ex}{$\tau_1=$}
\begin{tiny}
\raisebox{1ex}{ \pstree[levelsep=10ex]{\Tcircle{1} }{
   \Tcircle{1}

}}
\end{tiny}
\raisebox{-2ex}{\hspace{.1in} $\leftrightsquigarrow_1$ \hspace{.1in}
$\tau_2=$}
\begin{tiny}
\raisebox{1ex}{ \pstree[levelsep=10ex]{\Tcircle{1} }{
   \Tcircle{0}

}}
\end{tiny}
\raisebox{-2ex}{\hspace{.1in} $\leftrightsquigarrow_1$ \hspace{.1in}
$\tau_3=$}
\begin{tiny}
\raisebox{1ex}{ \pstree[levelsep=10ex]{\Tcircle{0} }{
   \Tcircle{0}}}
\end{tiny}
\end{center}
with $\tau_1=\tau(T')$ and $\tau_3=\tau(T'')$.  Using the note after
Definition \ref{definition:sim}, these trees correspond to the
sequence $T_1= T', T_2,$ and $T_3=T''$ of tableaux with
\begin{center}
\raisebox{3ex}{ $T_2=$ \;}
\begin{tiny}
\begin{tableau}
:.0>1>2\\
:>3>4\\
:>5>6\\
:^7^8\\
\end{tableau}
\end{tiny}
\raisebox{3ex}{.}
\end{center}
To show that $T_1 \leftrightsquigarrow' T_2 \leftrightsquigarrow'
T_3$, we need to find tableaux $S_1, S'_1, S_2,$ and $S'_2$ which
satisfy the equalities $MMT(S_1, T_1)=MMT(S'_1,T_2)$ and $MMT(S_2,
T_2)=MMT(S'_2,T_3).$ This is accomplished in our proof of Claim
\ref{claim} by using Proposition \ref{proposition:cyclestructure} to
construct tableaux with the required cycle structure:
\begin{center}
\raisebox{3ex}{$S_1=$ \;}
\begin{tiny}
\begin{tableau}
:.0>1>2\\
:^3^4^5^6\\
:;;;;\\
:^7>8\\
\end{tableau}
\end{tiny}
\raisebox{3ex}{ $S'_1=$ \;}
\begin{tiny}
\begin{tableau}
:.0>1>2\\
:^3^4>5\\
:;;>8\\
:^6^7\\
\end{tableau}
\end{tiny}
\raisebox{3ex}{$S_2=$ \;}
\begin{tiny}
\begin{tableau}
:.0>1>2\\
:^3^4^5^6\\
:;;;;\\
:^7^8\\
\end{tableau}
\end{tiny}
\end{center}
\vspace{.1in} and $S'_2=T''$.   We check that

$$
\raisebox{3ex}{$MMT(S_1, T_1)=MMT(S'_1,T_2)= $\;}
\begin{tiny}
\begin{tableau}
:.0.0>1>2\\
:.0>3>4\\
:^5>6\\
:;>8\\
:^7\\
\end{tableau}
\end{tiny}
$$

$$
\raisebox{3ex}{$MMT(S_2, T_2)=MMT(S'_2,T_3)=$ \;}
\begin{tiny}
\begin{tableau}
:.0.0>1>2\\
:.0>3>4\\
:^5>6\\
:;^8\\
:^7\\
\end{tableau}
\end{tiny}
$$

\end{example}

\end{document}